\documentclass[12pt]{article}
\usepackage{times}
\usepackage{charter}
\usepackage[charter]{mathdesign}
\usepackage{graphicx}
\usepackage{float}
\usepackage[font=rm]{caption}  
\usepackage{tabularx}
\usepackage{etoolbox}
\usepackage{array}
\usepackage{setspace}
\usepackage{amsmath}
\usepackage{amsfonts}
\usepackage{algorithm}
\usepackage[numbers]{natbib}
\setcitestyle{square}
\usepackage[margin=1.1in]{geometry}
\usepackage{setspace}
\AtBeginEnvironment{tabular}{\footnotesize}
\usepackage[affil-it]{authblk}
\date{}
\begin{document}
\title{\vspace{-1cm}{$N$-Dimensional curve following for solving \\ numerically systems of nonlinear equations}}
\author{\large Katerina G. Hadjifotinou}
\affil{\footnotesize Department of Mathematics, Aristotle University of Thessaloniki, Greece\\email: khad@math.auth.gr}
\maketitle
\vspace{-1cm}
\begin{abstract}
\noindent This paper presents a methodology for finding numerically, by means of curve-following, all real solutions of a general system of $n$ nonlinear equations in $n$ unknowns, within a given $n$-dimensional box. The main idea behind our method is a) to locate all parts of the curves formed by a selected subset of $n-1$ equations of the initial system, b) follow these parts numerically within the given $n$-box and c)  during this process, find all their intersection points with the hypersurface that represents the left-out equation of the initial system. With proper handling techniques, both stages (a) and (b) can be done with safety even when using a rapidly - but locally - convergent method such as Newton's method. Stage (c) on the other hand is theoretically straightforward and can be implemented by examining sign change and using bisection. However, improvement of performance with automatic step-size adaptation is also feasible. 

Since the choice of the left-out equation as well as the order of the unknowns, in some problems may affect the algorithm's performance, a secondary procedure for problem optimization has also been developed. This integrated algorithm has been successfully tested in all classical test-problems that are used in the literature (in total about 130 problems with dimension up to $n=10$) and achieved to locate approximations to all the known real solutions of each test-problem, with a user-provided accuracy. It has also proved to be more than 10 times quicker than Kuiken's \cite{Kuiken} 2D curve-following method, when compared in two-dimensional problems. 
\vspace{3.5\baselineskip}
\end{abstract}

\vspace{-1cm}


\section{Introduction}
\label{intro} Although there exists in literature a vast variety of numerical methods for finding one solution of a system of nonlinear equations, the problem of locating all real solutions in a given $n$-box has been addressed by a few authors in the past, and lately with renewed interest (see e.g. \cite{Chueca}).  Kearfott \cite{Kear87} provided a thorough analysis of the methodologies that could be used for such a problem and compared various interval, bisection and continuation methods, while at the same time investigating and presenting the restrictions and constraints of each. However he did not consider the idea of locating all the real solutions by $n$-dimensional curve following. This problem has only been studied in two dimensions by Kuiken \cite{Kuiken}, who presented a method of solution finding by curve following that is based on numerical integration. 

The basic notion behind Kuiken's method is that, the solutions of the system
\begin{equation*}
f(x,y)=0,~~g(x,y)=0
\end{equation*}
are merely the singular points of the system of ODEs  
\begin{equation*}
\frac{dx}{dt}=f(x,y),~~\frac{dy}{dt}=g(x,y).
\end{equation*}
However, Kuiken's method of curve-following is computationally expensive (as can be seen by comparison with our method in Sect.~6) and cannot be easily extended in more dimensions due to the rotations of the coordinate system that must take place at each step of the curve-following.

Despite the restrictions of Kuiken's method, the idea of extending the curve-following in many dimensions with the purpose of finding all the real solutions of a system of $n$ equations with $n$ unknowns can be in principle applicable. 

Let us assume that we need to calculate all real solutions of the $n$-dimensional system of nonlinear equations 
\begin{equation}
f_i(x_1,x_2,...,x_n)=0,~~i=1,..,n
\end{equation}
in the $n$-box defined by the known bounds $a_i, b_i$ such that 
$x_i\in[a_i,b_i], i=1,..,n$. 

In fact, one can always detach one equation from the system (1) and try to follow (in the $n$-dimensional space) the one-parametric curve formed by the remaining $n-1$ equations of (1) (provided the constant rank theorem holds) until one meets the hypersurface formed by the detached equation.

In this context, the main problems to be solved are a) how to locate all parts of the curves formed by the selected subset of $n-1$ equations and b) how to inexpensively follow each of these parts within the given $n$-box. Apart from the above obvious issues, another question that arises from this specific approach is, how to select which equation to leave out from the original system and consequently, how the choice of the subset of equations to follow, affects the overall algorithmic performance.

In the sections that follow, we address all the above problems, together with that of the implementation of the final stage c) of efficiently locating, during the curve following, the intersection points of the followed curve with the hypersurface of the left-out equation. In specific, in Sect.~2 we explain stage (a) and in Sect.~3 we explain stage (b) together with stage (c) of our method, since these two stages are executed in combination. In Sect.~4, the procedure is illustrated through two examples. In Sect.~5, a thorough discussion is made on the dependence of our method's performance on the ordering of the system's rows and variables, and an auxiliary procedure to improve this performance through changing of ordering is presented. Finally, in Sect.~6, applications of this integrated method to the most representative of the known test-problems are presented and a comparison of our algorithm's 2-dimensional case with the method of Kuiken \cite{Kuiken} is elaborated. Suggestions for further extensions and possible improvements of the method are discussed in Sect.~7.

\section{Algorithm's stage (a): Locating all parts of a curve in $\mathbb{R}^n$}
\label{sec2} 
The first question that arises in $n$-dimensional curve following, is obviously, how to locate all the parts of the curve that need to be followed. An $n$- dimensional exhaustive grid search is computationally prohibitive due to its large algorithmic complexity.
A feasible approach however, is to choose one variable (e.g. $x_n$), consider it as "running" variable (that is, as the curve's parameter) and examine the intersections of the curve with "slices" at this dimension (i.e. $x_n$-slices). This procedure was also followed by Kuiken \cite{Kuiken} in his two-dimensional approach. The algorithm we have created for this purpose is presented below:\\
{\vskip 0.1cm}
\noindent{\bf Algorithm 2.1: LocateCurveParts}\\{\it\  
Step 1: Input \{F, J, n, lower, upper, stepx, stepz, acc$_1$, acc$_2$, step, thresh\}\\
Step 2: ListOfCurveParts=\{\}\\
Step 3: $F_l$=F(1:n-1); $~~F_u$=F(n); $~~~$ $J_l$=J(1:n-1;1:n-1)\\
\hspace*{37pt} $lower_l$=lower(1:n-1); $~$ upper$_l=$upper(1:n-1)\\
Step 4: npoints=floor((upper$_l-lower_l$)/stepx+1);~meshpoints=prod(npoints)\\ 
\hspace*{37pt} A=Rmesh(npoints, lower$_l$, stepx)\\	
Step 5: z$_0$=lower(n)\\ 
Step 6: Repeat Steps 7 - 19 while z$_0$<=upper(n)\\ 
\hspace*{18pt} Step 7: $F_{l0}=F_l(x_n=z_0); ~~~ J_{l0}=J_l(x_n=z_0)$\\
\hspace*{18pt} Step 8: CurveStartingPoints=$[]$\\
\hspace*{18pt} Step 9: Repeat Steps 10 - 12 for i=1:meshpoints\\
\hspace*{38pt} Step 10: $v_0=A(i,:)$\\
\hspace*{38pt} Step 11: $[v, found]=Newton(n-1,F_{l0},J_{l0},v_0,acc_1$)\\
\hspace*{39pt} Step 12: If found: CurveStartingPoints=Append(CurveStartingPoints,v,$acc_1$)\\
\hspace*{38pt} Step 13: Repeat Steps 14 - 18 for each v in CurveStartingPoints\\
\hspace*{52pt} Step 14: $vz_0=[v~z_0]$\\
\hspace*{52pt} Step 15: exists=Belongs($vz_0$,~ListOfCurveParts)\\
\hspace*{52pt} Step 16: If not(exists): Do Steps 17 - 18\\ 
\hspace*{72pt} Step 17: $[$ListOfCurveParts, ListOfSolutions$]$=FollowCurve($F_l, F_u, J_l$,\\
\hspace*{84pt} n, lower, upper, step, $vz_0, acc_1, acc_2,$ thresh, ListOfCurveParts,$~$ListOfSolutions) \\
\hspace*{72pt} Step 18: $[$ListOfCurveParts, ListOfSolutions$]$=FollowCurve($F_l, F_u, J_l$,\\
\hspace*{84pt} n, lower, upper, -step, $vz_0, acc_1, acc_2,$ thresh, ListOfCurveParts,$~$ListOfSolutions) \\
\hspace*{18pt} Step 19: $z_0=z_0+$stepz\\
Step 20:  Output \{ListOfSolutions\}}
{\vskip 0.4cm}
\noindent Each step of the algorithm is thorougly explained below.

{\it Step 1:} $F$ is the vector of functions in variables $x_1,..,x_n$ of system (1), $J$ is the Jacobian matrix of $F$, $lower$ is a vector of all the lower bounds $a_i$, and $upper$ is a vector of all the upper bounds $b_i$ that define the $n$-box into which we are searching for the solutions of (1). $stepx$ is the step interval-length for the creation of subintervals of $[a_i, b_i], i=1,..,n-1$ (see $Step~4$ below). Before executing Algorithm 2.1 and without loss of generality, we assume that our equations and variables are already ordered in such a way that,  the last equation of (1) will be left out of the system and the last variable ($x_n$) will be "sliced".  $stepz$ is the step interval-length for the creation of subintervals of $[a_n, b_n]$, and therefore it is the step-size of this "slicing". The optimal way of variable and row ordering, and a routine to do this automatically, before the call of Algorithm 2.1, will be extensively discussed in Sect.~5. Finally, $acc_1$ is the user-specified accuracy with which curve points will be found and $acc_2$ is the accuracy of solution finding. $step$ is the initial step used for curve-following and is passed on to algorithm $FollowCurve$, together with the input parameters $acc_2$ and $thresh$. The latter will be discussed in the next section.

{\it Step 2:} For reasons of efficiency, the algorithm stores the parts of the curves followed so far in a list, so that these curves are not visited again. Initially this list is empty.

{\it Step 3:} Split the vector of functions $F$ into two parts: $F_l$, the functions of the first \mbox{$n-1$} equations of (1) and $F_u$, the last equation of (1). Correspondingly, denote with $J_l$ the submatrix of the Jacobian $J$ that corresponds to $F_l$. Also denote with $lower_l, upper_l$ the vectors with the first $n-1$ elements of the bounds' vectors $lower$ and $upper$ respecively.

{\it Step 4:} Create an $(n-1)$-dimensional mesh of starting points with step $stepx$ within the $(n-1)$-box defined by the vectors $lower_l, upper_l$. The number of these points is $meshpoints$. These points are stored in a 2-dimensional matrix A, whose number of rows is $meshpoints$ and each row contains the $n-1$ coordinates of a starting point. The routine that creates the starting points' matrix A is especially created for this purpose and is called $Rmesh$.

{\it Step 5:} Take the first "slice" (called $z_0$) of the last variable $x_n$, whose value is the lower bound $a_n$ of $[a_n,b_n]=[${\it lower(n), upper(n)}$]$.

{\it Step 6:} Repeat Steps 7 - 19 for each "slice" $z_0$ of $x_n$, until we reach its upper bound $b_n$.

{\it Step 7:} Substitute $z_0$ to $x_n$ into the system $F_l=0$ and create the $(n-1)$-dimensional subsystem $F_{l0}=0$. Substitute $z_0$ to $x_n$ also in the submatrix $J_l$ and create $J_{l0}$.

{\it Steps 8-12:} For each starting point $v_0$ (i.e. for each row) of matrix $A$, use this point as a starting point for any locally convergent method such as Newton's method, in order to find a solution of the $(n-1)\times(n-1)$ subsystem $F_{l0}=0$ with accuracy $acc_1$. If Newton's method converges to a solution $v$ (and so the flag-variable $found$ is returned as $true$), then this $(n-1)$-dimensional solution, when appended with the coordinate $x_n=z_0$, will be a point of the curve of the subsystem $F_l=0$ of (1) for the current $x_n$-slice $z_0$. We add the solution vector $v$ to the matrix $CurveStartingPoints$, which is a 2-dimensional matrix that is initially empty and at the end of the loop of {\it Step~9} will contain at most as many points as matrix $A$ (the solutions of the subsystem $F_{l0}=0$ that correspond to the mesh starting points in $A$). If Newton's method diverges (and so the flag-variable $found$ is returned as $false$), no point is appended to the matrix $CurveStartingPoints$. The inclusion of the solution $v$ to the matrix $CurveStartingPoints$ is done with the routine $Append$, which also checks if the certain solution coincides (with accuracy $acc_1$) with any already existing in $CurveStartingPoints$ solution, to which Newton's method may have converged when starting from a previously used starting point $v_0$ of $A$. By the end of the loop of {\it Step~9}, the number of rows of $CurveStartingPoints$ must be the number of all the different curve-parts of the subsystem $F_l=0$ that exist at the "slice" $x_n=z_0$.

{\it Steps 13-18:} Use each curve-point found in steps 9-12 (after adding to it the last coordinate $x_n=z_0$) as a starting point (named $vz_0$) to follow this curve upwards and downwards into the user-given $n$-box. Store the curve-points found, updating the $ListOfCurveParts$. However, before starting the curve following, we need to check if $vz_0$ already exists in a curve part found in a previous iteration of the loop of $Step~ 13$. This task is accomplished by the routine $Belongs$, which returns the flag-variable $exists$. If $exists$ is $true$, then nothing needs to be done, since this curve part has already been followed. In this way we avoid redundant curve-following. Finally, the algorithm $FollowCurve$ starts from the curve-point $vz_0$ and follows the curve upwards ($Step 17$) and then downwards ($Step 18$), with step-size $step$, simultaneously appending the $ListOfCurveParts$ and investigating the existence of intersection points with the surface of the last equation $F_u=0$. $FollowCurve$ will be presented in the next section.

{\it Step 19:} Increase the $z_0$-slice by $stepz$ and proceed to the next iteration of the loop of $Step~6$. At each $z_0$-slice, new curve parts may be added to the $ListOfCurveParts$. At the end of this loop, the list $ListOfCurveParts$ should contain all the parts of the curves of the subsystem $F_l=0$ that belong to the wanted $n$-box.

{\it Step 20:} The algorithm returns the $ListOfSolutions$, that is updated each time the routine $FollowCurve$ is executed. With proper selection of step-sizes, by the end of the loop of $Step~6$, the  $ListOfSolutions$ should contain all the solutions of (1) in the given $n$-box.

The routines $Rmesh$, $Append$ and $Belongs$ that are called by Algorithm 2.1, are auxiliary routines especially constructed for use by our method and will not be presented here since they contain elementary calculations, but can be provided by the author upon request. On the contrary, the algorithm $FollowCurve$ comprises an important part of the method and therefore it will be explained in detail in the next section.

\section{Algorithm's stage (b): Following a curve in $\mathbb{R}^n$}
\label{sec3} 

In the present section we present the algorithm $FollowCurve$, which starts from a curve-point of the list $CurveStartingPoints$ that was created by Step 12 of Algorithm 2.1 and follows this curve in $\mathbb{R}^n$ for increasing or decreasing $x_n$-values, depending on the sign of the input parameter $step$. During the curve-following process, an investigation of possible intersection points with the surface of the last equation $F_u=0$ is carried out.\\
{\vskip 0.1cm}
\noindent{\bf Algorithm 3.1: FollowCurve}\\{\it\
Step 1: Input \{$F_l, F_u, J_l$, n, lower, upper, step, $vz_0, acc_1, acc_2,$ thresh,$~$ListOfCurveParts,\\\hspace*{70pt}$~$ListOfSolutions\}\\
Step 2: Repeat Steps 3 - 16 while lower(n)$\leq  vz_0(n) \leq$ upper(n)\\
\hspace*{20pt} $Step 3: ListOfCurveParts=Append(ListOfCurveParts,vz_0,acc_1$)\\
\hspace*{17pt} Step 4: $u_0=F_u(vz_0$)\\
\hspace*{17pt} Step 5: If $|u_0|\leq acc_2$: Do Steps 6 - 8\\
\hspace*{38pt} Step 6: ListOfSolutions=Append(ListOfSolutions,$vz_0,acc_1$)\\
\hspace*{38pt} Step 7: $[vz_0, done]$=ProceedOneStep($F_l, F_u, J_l, n, step, vz_0, acc_1$,thresh)\\
\hspace*{38pt} Step 8: If not(done): Exit Loop of Step 2\\
\hspace*{18pt} Step 9: Else Do Steps 10 - 16\\
\hspace*{38pt} Step 10: $[vz_1$, done$]$=ProceedOneStep($F_l, F_u, J_l$, n, step, $vz_0, acc_1$,thresh)\\
\hspace*{38pt} Step 11: If not(done): Exit Loop of Step 2\\
\hspace*{38pt} Step 12: $u_1=F_u(vz_1)$\\
\hspace*{38pt} Step 13: If $|u_1|>acc_2$ and $u_0\cdot u_1<0$: Do Steps 14 - 15\\
\hspace*{52pt} Step 14: [solution found]=Bisection($F_l, F_u, J_l, n, vz_0, vz_1, acc_2$)\\
\hspace*{52pt} Step 15: If found Add solution to ListOfSolutions\\
\hspace*{38pt} Step 16: $u_0=u_1$;$~~$ $vz_0=vz_1$\\
Step 17:  Output \{ListOfCurveParts, ListOfSolutions\}
}
{\vskip 0.3cm}
\noindent Each step of the algorithm is thorougly explained below.

{\it Step 1:} $F_l, F_u, J_l, n, lower$ and $upper$ are as explained in the previous section.
$step$ is the initial $x_n$-step used for curve following, with $x_n$ being the running variable. This step might decrease in order to meet the proximity demands of two successive curve points - see description of algorithm $ProceedOneStep$ for more details. $vz_0$ is the curve's starting point. $acc_1$ is the wanted accuracy of curve-point approximation. It must be relatively high (e.g. $10^{-10}$) for the curve to be followed correctly. On the contrary, $acc_2$ (the wanted accuracy of the approximation of solutions used by $Bisection$ routine) needs not be high (e.g. $10^{-4}$), since our method aims to provide initial approximations only of all the solutions of (1) within a given $n$-box. These approximations can subsequently be used as starting values by any high-precision single-solution finding algorithm. $thresh$ is a threshold step-size value that is passed on to algorithm $ProceedOneStep$ and will be discussed later in this section. Finally, both $ListOfCurveParts$ and $ListOfSolutions$ are updated at each call of $FollowCurve$, and that is why they are both Input and Output to the algorithm.

{\it Step 2:} All the steps that follow must be repeated as long as the running variable $vz_0(n)=x_n$ stays within $[a_n, b_n]=[${\it lower(n), upper(n)}$]$.

{\it Step 3:} The starting curve point $vz_0$ was calculated by Newton method (see steps 11 and 14 of Algorithm 2.1) and satisfies the equations $F_l=0$ for a certain "slice" $x_n=z_0$. Since it belongs to the curve, it must be added in the $ListOfCurveParts$. After each new loop iteration, $vz_0$ will hold the new curve point (see steps 7, 10 and 16 above) and therefore it must be also added in the $ListOfCurveParts$.

{\it Step 4:} In order to investigate the existence of solutions of (1), we calculate the value $u_0$ of its last equation $F_u=0$ at each new curve point $vz_0$.

{\it Steps 5-8:} We take into account the rare case that in the current curve-following step we have stepped directly on a solution (actually, with proximity $acc_2$). In that case, the approximate solution $vz_0$ is added to the $ListOf$ $Solutions$ and we procceed to the next step of curve following, calling the routine $ProceedOneStep$ (see Algorithm 3.2 below) and calculating a new $vz_0$. In this way, if the step-size $step$ is sufficiently small, we can also locate clustered solutions. If the routine $ProceedOneStep$ does not manage to proceed further, this means that the curve part we are following is not continued in the same $x_n$ direction (either increasing or decreasing, depending on the sign of $step$) and algorithm $FollowCurve$ stops.   

{\it Steps 9-16:} If we have not already found a solution, we have to compare the value $u_0$ of $F_u$ at the last calculated $vz_0$ with the value $u_1$ at the point $vz_1$ which is the new point calculated by calling $ProceedOneStep$. Of course, as we have mentioned before, if 
$ProceedOneStep$ fails to calculate $vz_1$, the algorithm stops. If we do have a $u_1$ however, there is again the chance that this is absolutely small enough to consider that we have approximately found a solution ($vz_1$). In that case, we rename the solution point to $vz_0$ and the corresponding $F_u$ value to $u_0$ ($Step~16$) and proceed the loop iterations, so that steps 5-8 are executed, in order to handle this new solution and then to move on.

If $u_1$ is not small enough for $vz_1$ to be considered a solution, we check if it is eterosign of $u_0$ ($Step~13$). Then there may exist a solution of (1) with $x_n$ between $vz_0(n)$ and $vz_1(n)$. This solution is found using a bisection-type algorithm ($Step~14$) that will be presented later in this section. If the routine $Bisection$ achieves to locate the solution, then the latter is added to the $ListOfSolutions$ ($Step~15$) and the algorithm proceeds after setting $vz_1$ to $vz_0$ and $u_1$ to $u_0$ ($Step~16$). 

On the other hand, if $u_1$ is of the same sign as $u_0$, we proceed directly to $Step~16$. A discussion can be done at this point concerning the possibilities of increasing or decreasing the step-size $step$ of curve-following according to whether $u_1$ increases or decreases in absolute value, compared to $u_0$ (that is, according to whether we are distancing from or approaching to a solution). However, as also mentioned in Sect.~7, tests that we have performed for this purpose, led us to the conclusion that this alteration of $step$ affects the adaptation of $h$ which is done anyway for the needs of precise curve-following in algorithm $ProceedOneStep$ presented below. As a result, the overall performance of the algorithm is not improved.

In what follows, we present the algorithn $ProceedOneStep$.

{\vskip 0.3cm}
\noindent{\bf Algorithm 3.2: ProceedOneStep}\\{\it\
Step 1: Input \{$F_l, F_u, J_l, n,~step,~vz_0,~acc_1,~thresh$\}\\
Step 2: $v_0=vz_0($1:n-1$)$; $~z_0=vz_0(n); \quad$ h=step\\
Step 3: $z_0=z_0+h; \quad F_{l0}=F_l(x_n=z_0); \quad J_{l0}=J_l(x_n=z_0)$\\
Step 4: $[$v, found$]$=Newton(n-1,$F_{l0},J_{l0},v_0,acc_1); ~~ vz=[v z_0]$\\
Step 5: NotAcceptable=not(found) or (found and norm(v-v$_0$)>step)\\
Step 6: Repeat Steps 7 - 10 while NotAcceptable and $h\geq thresh$ \\
\hspace*{18pt} Step 7: h=h/2;  $~ z_0=z_0+h; ~~ F_{l0}=F_l(x_n=z_0); \quad J_{l0}=J_l(x_n=z_0)$\\
\hspace*{18pt} Step 8: $[$v, found$]$=Newton$(n-1,F_{l0},J_{l0},v_0,acc_1)$\\
\hspace*{18pt} Step 9: NotAcceptable=not(found) or (found and norm(v-v$_0$)>step)\\
\hspace*{18pt} Step 10: done=not(NotAcceptable); $~$vz=$[$v$~z_0]$\\ 
Step 11: Output \{vz, done\}
}	
{\vskip 0.3cm}
\noindent In brief, what Algorithm 3.3 does, is to advance curve-following from $vz_0$ to $vz$. If Newton method diverges when starting from $v_0=vz_0(1:n-1)$, or even if it returns a solution $v$ of the $(n-1)\times (n-1)$ subsystem of (1) that is far away from the initial point $v_0$, then we consider that curve-following is not done correctly and $v$ is not acceptable ($Step~5$). In that case, we enter a loop (steps 7-10) of repeating the procedure while halving the step-size $h$, until either we find an acceptable solution $v$ or $h$ becomes so small that we are practically not advancing the curve-following. The lower bound of $h$ is controlled by the user-defined threshold parameter $thresh$, that has been passed to $ProceedOneStep$ through the calls of $LocateCurveParts$ and $FollowCurve$.
After the end of this loop, the algorithm either returns an acceptable next curve-point $vz$, or it outputs $false$ through the flag variable $done$.

As a last structural element of our method, we present the algorithn $Bisection$ that is called at $Step~14$ of algorithm $FollowCurve$. This algorithm applies the principle of bisection in the $n$-th dimension of our problem and uses Newton method for calculating the $n-1$ coordinates of the point with middle $x_n$.
 
{\vskip 0.3cm}
\noindent{\bf Algorithm 3.3: Bisection}\\{\it\
Step 1: Input \{$F_l, F_u, J_l, n, vz_0, vz_1, acc_1, acc_2$\}\\
Step 2:~zm=($vz_0(n)+vz_1(n))/2;F_{l0}=F_l(x_n=zm);J_{l0}=J_l(x_n=zm);v_0=vz_0$(1:n-1)\\
Step 3: $[$midpoint, found$]$=Newton(n-1,$F_{l0},J_{l0},v_0,acc_1)$\\
Step 4: If not(found): End\\
Step 5: vm=$[$midpoint zm$]; ~~gm=F_u(vm); ~~ga=F_u(vz_0); ~~gb=F_u(vz_1)$\\
Step 6: Repeat Steps 7 - 19 while $|gm|>acc_2$ and found\\
\hspace*{18pt} Step 7: If $ga\cdot gm<0$ and $|gm|\leq$|gb|: Do Steps 8-11\\
\hspace*{35pt} Step 8: $vz_1$=vm$; ~~zm=(vz_0(n)+vz_1(n))/2$\\
\hspace*{35pt} Step 9: $F_{l0}=F_l(x_n=zm); ~~~ J_{l0}=J_l(x_n=zm)$\\ 
\hspace*{35pt} Step 10: [midpoint, found]=Newton(n-1,$F_{l0},J_{l0},v_0,acc_1)$\\
\hspace*{35pt} Step 11: If not(found): End.\\
\hspace*{35pt} Step 12: vm=[midpoint zm]; $~~gb=gm; ~~gm=F_u$(vm)\\
\hspace*{18pt} Step 13: Else If gm$\cdot gb<~0$ and $|gm|\leq$|ga|: Do Steps 14-18\\
\hspace*{35pt} Step 14: $vz_0=vm; ~~zm=(vz_0(n)+vz_1(n))/2$\\
\hspace*{35pt} Step 15: $F_{l0}=F_l(x_n=zm); ~~~ J_{l0}=J_l(x_n=zm)$\\ 
\hspace*{35pt} Step 16: [midpoint, found]=Newton(n-1,$F_{l0},J_{l0},v_0,acc_1)$\\
\hspace*{35pt} Step 17: If not(found): End\\
\hspace*{35pt} Step 18: vm=$[$midpoint zm$]; ~~ga=gm; ~~gm=F_u$(vm)\\
\hspace*{18pt} Step 19: Else: found=False; Exit Loop of Step 6 \\
Step 20: If |gm|$\leq acc_2$: solution=vm\\
Step 21: Output \{solution, found\}
}\\

\noindent Each step of the algorithm is explained below.

{\it Step 1:} Apart from the obvious $F_l, F_u, J_l$ and $n$ that are passed to the routine through $FollowCurve$, the important input here is the vectors $vz_0$ and $vz_1$ that are the two curve points between which a possible intersection point with the surface $F_u=0$ exists, since $F_u(vz_0)$ and $F_u(vz_1)$ are eterosign (see steps 13-14 of Algorithm 3.1). As we have mentioned before, $acc_1$ is the accuracy of curve-following, while $acc_2$ is the accuracy of solution finding.

{\it Step 2:} We locate the middle point $zm$ that concerns the last dimension ($x_n$) of the two vectors $vz_0$ and $vz_1$ and substitute this point to the functions $F_l$ and $J_l$.

{\it Step 3:} By calling the Newton routine to solve this system for the middle "slice" $zm$ and for the starting point $v_0=vz_0(1:n-1)$, we find the curve-point that belongs to this middle $x_n$-slice. 

{\it Step 4:} If Newton fails, given the fact that it has already managed to locate correctly the point $v_1=vz_1(1:n-1)$ at $Step~10$ of $FollowCurve$, this indicates that the curve possesses a singularity instead of an intersection point, between the $x_n$-slices $z_0=vz_0(n)$ and $z_1=vz_1(n)$. As a result, the algorithm stops returning $false$.

{\it Steps 5-20:} The $x_n$-middle point of the curve has been found and we calculate the value $gm$ of $F_u$ at this point. If $gm$ is absolutely small enough to meet the user-specified $acc_2$ demand, then the solution is considered found and it is returned by the algorithm. Otherwise, the procedure must be repeated for the left (steps 7-12) or right (steps 13-18)  subinterval of $[z_0, z_1]$, depending on where the change of sign happens. The iterative procedure stops, either if, by the successive halving of the $x_n$-interval width, we finally locate the intersection point, or if Newton method stops converging, or if we reach a point whose $F_u$ absolute value increases instead of decreasing. In the two last cases, the algorithm returns $false$.    

\section{Examples}
\label{sec4} 

The exact way Algorithms 2.1-3.3 work, can be better illustrated through the following examples. Let us first consider the $3\times 3$ system, resulting from the $n=3$ version of the Trigonometric function problem of Spedicato \cite{Spedicato}, also presented as Problem no. 26 in \cite{MGH}:
\begin{equation}
f_i(x_1,x_2,x_3)=3-\sum_{j=1}^{3}\cos{x_j}+i(1-\cos{x_i})-\sin{x_i}=0, ~~i=1,..,3.\\
\end{equation}
Let us assume that we want to locate all the real solutions of (2) within the three-dimensional box $[-10,10]^3$. As seen in Figure~1, each of the three equations represents a surface that consists of ellipsoids in a mesh arrangement. In the denoted box, there are 27 ellipsoids in each surface, which, when combined altogether, form intersecting clusters of three, as seen in a detail for the smaller $[-2,2]^3$ box in Figure~2. Each triple of intersecting ellipsoids leads to the existence of two intersection points. As a result, the system (2) possesses 2 solutions in the box $[-2,2]^3$ and 54 solutions in the box $[-10,10]^3$.

\begin{figure}[H]
\centering
\includegraphics[width=1\textwidth]{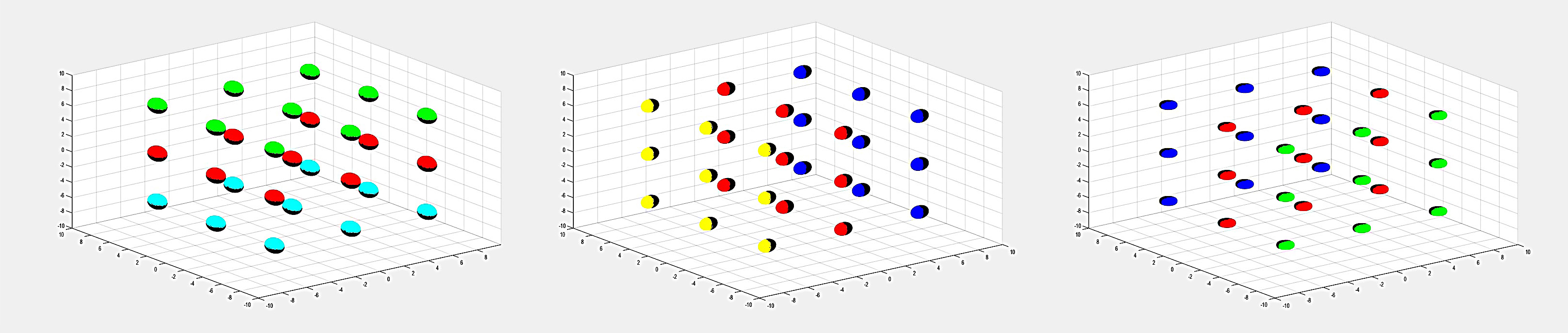}
\caption{The surfaces that correspond to each one of the equations (2) in the box $[-10,10]^3$.}
\label{fig1}
\end{figure}
\begin{figure}[H]
\centering
\includegraphics[width=0.5\textwidth]{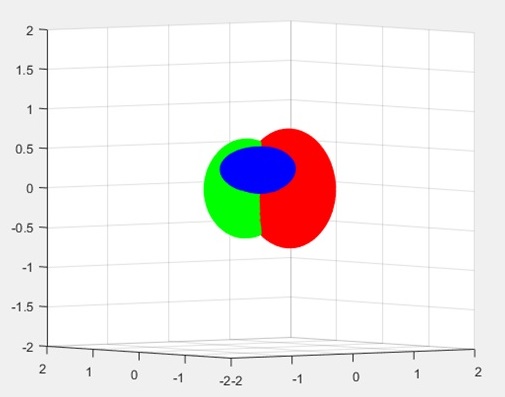}
\caption{All three surfaces of (2) ($f_1$:red, $f_2$:green, $f_3$:blue) in the sub-box $[-2,2]^3$.}
\end{figure}

Let us execute Algorithm 2.1 with arguments $F=(f_1, f_2, f_3)$ of (2), $J$ the Jacobian matrix of $F$, $n=3$, $lower=(-10, -10, -10),~upper=(10, 10, 10)$, {\it stepx}=1, {\it stepz}=1, $acc_1=10^{-10},~acc_2=10^{-4}$, $step=0.1$ and $thresh=0.1$. With $stepx=1$, the mesh points that are created as starting values of Newton method in the 2-dimensional sub-box $[-10, 10]^2$ of $[-10, 10]^3$, have coordinates that combine the 21 points \{-10, -9,.., 9, 10\} in the $x_1$ direction with the same points in the $x_2$ direction. Thus, a mesh of 441 initial points is created. 

As concerns the "running" variable $x_3$, since we have also chosen $stepz=1$, 21 $x_3$-slices are formed with $z_0\in \{-10, -9,.., 9, 10\}$. For each slice, Algorithm 2.1 attempts to find curve points of the subsystem of the first two equations of (2). E.g. for the first slice $x_3=z_0=-10$, as can also be verified by Figure~1, no real solutions exist for the $2\times 2$ subsystem of (2)
\begin{equation}
\begin{aligned}
f_1(x_1,x_2)=4-2\cos{x_1}-\cos{x_2}-\sin{x_1}-\cos{z_0}=0\\
f_2(x_1,x_2)=5-\cos{x_1}-3\cos{x_2}-\sin{x_2}-\cos{z_0}=0\\
\end{aligned}
\end{equation}
and therefore, no curve points are found. As proceeding the iterations of the loop of $Step~6$ of Algorithm 2.1, the first "slice" for which curve points are located  is slice $z_0=-6$. For this $z_0$ value, when using the 441 mesh starting points that we have mentioned before, the application of Newton method to system (3) results to the location of 18 different real solutions of (3) in the 2-box $[-10, 10]^2$. This is correct since, the "slice" $z_0=-6$ meets only 9 triples of the intersecting ellipsoids of (2) (see Figure~1) and each triple contains one closed elliptical curve as intersection of the surfaces $f_1$ and $f_2$ (the red and green surfaces in Figure~2 respectively). Therefore a $x_3$-slice will meet each elliptical curve at maximum two points. 

Before proceeding to the next $z_0$-slice, the routine $FollowCurve$ is called twice, once for increasing $x_3$ with $step=0.1$ ($Step~17$ of Algorithm 2.1), and once for decreasing $x_3$ with $step=-0.1$ ($Step~18$ of the same algorithm), and this is done for each one of the 18 located curve points. As a curve is followed upwards or downwards, investigation of intersection points with the surface $f_3(x_1, x_2, x_3)=0$ is simultaneously run (Algorithm 3.1). In our example, for $z_0=-6$, 9 of the 18 curve points produce curve-parts that have 2 intersection points with the surface of the third equation, while the other 9 curve points produce curve parts that have no intersection points with $f_3=0$. As a result, we obtain 18 solutions when starting from the 18 curve points that correspond to the $z_0=-6$ slice.

As the $z_0$-slices proceed after $z_0=-6$, we discover that there exist only two other slices, $z_0=0$ and $z_0=6$, for which the 441 mesh starting points obtained by $Step~4$ of Algorithm 2.1, make Newton produce solutions of (3). Again, as expected from Figure~1, 18 curve points are added by Algorithm 2.1 to the list $CurveStartingPoints$ for each of the two $z_0$-slices. Subsequently, during the $FollowCurve$ process, these two sets of 18 curve starting points each, lead to curve parts that have 18 intersection points with $f_3=0$ for the set with $z_0=0$ and 18 more intersection points for the set with $z_0=6$. As a result, our method outputs correctly all the 54 solutions of (2) that exist in the 3-box $[-10, 10]^3$. 

It is worth noting here, that some solutions are found while following a curve branch with ascending $x_3$ during the $FollowCurve$ procedure, while others are located while following a curve branch with descending $x_3$. From the 54 curve starting points that are produced in total for the three $z_0$-slices $z_0=-6$, $z_0=0$ and $z_0=6$, half of these points belong to curve parts that have no intersection points with $f_3=0$, and the other half belong to curve parts that have 2 intersection points each with $f_3=0$, as we have also mentioned before. 

In the above example, the geometry of the curves (2) is such, that each new starting point found by Algorithm~2.1 belongs to a different curve branch. As a result, the necessity of storing the visited curve points in the list $ListOfCurve$ $Parts$ in Algorithm~2.1 and the additional check with the routine $Belongs$, so as to avoid the need of refollowing a branch, is not made evident through problem (2). For this purpose, as well as for highlighting some other issues, let us also examin, in the 2-box $[-2, 2]^2$, the following 2D example:

\begin{equation}
\begin{aligned}
f_1(x_1,x_2)=\sin(x_1^2+2x_2^2)=0\\
f_2(x_1,x_2)=\tan(x_1^2-2x_2^2)=0\\
\end{aligned}
\end{equation}  

 Using as input for Algorithm~2.1 the parameters {\it stepx}=0.5, {\it stepz}=0.5, $step=0.1$ and $thresh=0.1$, we consider as "running" variable the variable $x_2$ and follow the curves of the first equation (black curves in Figure~3a) in order to meet the curves of the second equation (red curves in Figure~3a). For this $stepx$ value we acquire the 1D mesh of 9 starting points \{-2, -1.5, ..., 1.5, 2\} in the $x_1$ direction. Since for the "running" variable $x_2$  we have also chosen $stepz=0.5$, 9 $x_2$-slices are formed with $z_0\in \{-2, -1.5, ..., 1.5, 2\}$.

For each $z_0$-slice and using the 9 mesh points as starting points for the one-dimensional Newton method to solve for $x_1$ the single nonlinear equation $f_1(x_1, z_0)=0$, two different solutions are found by Newton - apart from the $z_0=0$ slice for which three solutions are found. These solutions (in total 19) are the $CurveStartingPoints$ used by algorithm $FollowCurve$ and are marked with blue stars in Figure~3a. Using each one of these curve starting points and following its curve with ascending and then with descending $x_2$, as shown in Figure~3b, all black curve parts are covered and consequently all 27 solutions of system (4) in the 2-box $[-2, 2]^2$ are found (together with the isolated solution (0,0)). 
For the needs of performing an accurate following process, the step-size is often automatically halved (see $Step~7$ of Algorithm 3.1). As a result, a clustering of stars (i.e. visited points) appears in Fig.~3b, at the zero-derivative parts of the black curve (around $x_1=0$).  

\begin{figure}[h]
	\centering
	\includegraphics[width=1\textwidth]{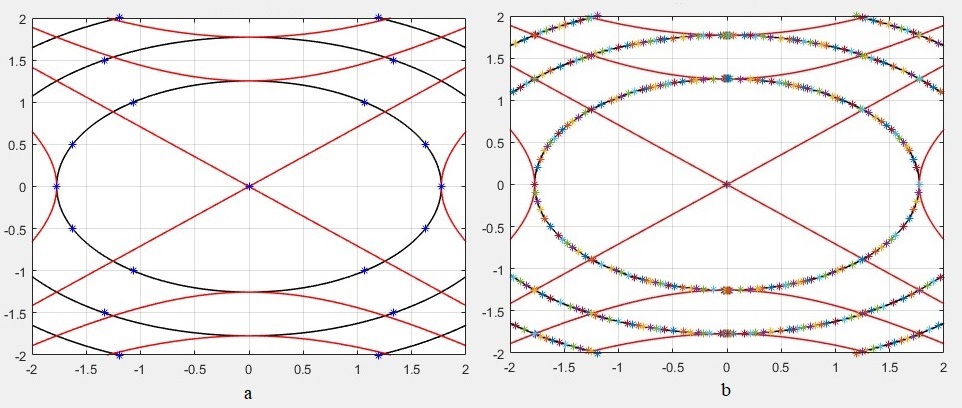}
	\caption{a. The curves of the equations (4) ($f_1$:black, $f_2$:red) and the starting points of each curve branch found by Algorithm 2.1 (blue stars). b. After the execution of Algorithm 3.1, all visited curve points of $f_1=0$ are marked with stars. The isolated point (0,0) is also visited. Notice also the clustering of points around the four difficult-to-follow zero-derivative curve parts at $x_1=0$ that demand step-size halving.}
	\label{fig2}
\end{figure}

As seen from Fig.~3a, many of the 19 curve starting points belong to the same curve parts. For instance, as regards the inner ellipse, from the 10 starting points that belong to it, two of them (one with negative and one with positive $x_1$) suffice as starting points to cover the whole ellipse, when we follow it with running $x_2$ ascending and then descending. Therefore, in order to avoid redundant curve-following, an obvious strategy is to store the points of the already visited curve parts and test if a new starting point belongs to these parts (routine $Belongs$ of Algorithm 2.1) before implementing the curve-following. This is the case for the starting points with $z_0$ from $-0.5$ up to $z_0=1$, since all of them - apart from the point $(0,0)$ - belong to the inner ellipse of Fig.~3a that is already followed, commencing from the two curve-starting points with $z_0=-1$. Each of the rest of the points leads to a number of solutions when ascending $x_2$ and/or when descending $x_2$, and three of them (the points (1.3354, -1.5), (1.0685, -1) and (1.3354, 1.5)) lead to a solution that is already approached from their symmetric point with negative $x_1$. Finally, the starting point $(0,0)$ is itself an isolated solution and comprises a visited curve-part on its own. 

\section{Changing the equation and variable ordering}
\label{sec5} 

Since the main idea of our curve-following methodology is based on excluding an equation from the system (1) and on choosing a particular dimension as the one that will be "sliced", this leads automatically to the conclusion that the choice of the left-out equation as well as the choice of the "running" variable plays an important role on the algorithm's performance.
In an attempt to optimize performance, firstly in term of the ability to locate all the solutions in a $n$-box when following a specific curve, and secondly to achieve this with the largest step-size possible - i.e. using the minimum possible number of "slices", we first implemented two more routines. The routine $SwapRows$ automatically changes the order of equations by swapping any row between row~1 and row~$n-1$ (selected by the user) with the last row, since it is always the last row that is left out by the algorithm. Of course, the Jacobian $J$ is also altered accordingly. The routine $SwapColumns$ automatically changes the order of variables in the input functions $F$ and $J$, on user demand, again by swapping any variable between $x_1$ and $x_{n-1}$ with the variable $x_n$, since $x_n$ is always considered as the "running" variable. These routines are not presented here since they are elementary, but can be provided by the author upon request.

After adding the routines $SwapRows$ and $SwapColumns$, we tried to execute our algorithm for each one of our test-problems, with various choices of the left-out equation and running variable, attempting to spot any relationship between the optimal selection of equation or variable and the properties of the Jacobian of $F$. The Jacobian can be of use, since it is a quantity that provides information on the type of dependence of each function $f_i$ in (1) with respect to each variable. For instance, in a $2\times 2$ system, if the bivariate function $f_1(x_1, x_2)$ is nonlinearly dependent on variable $x_1$ and linearly dependent on variable $x_2$, it will be much easier to use $x_1$ as the running variable instead of (the default) $x_2$, since, substituting any $z_0$ value of a $x_1$-slice to $f_1$ will lead to a simple linear equation: $f_1(z_0,x_2)=0$, to be numerically solved and followed in $\mathbb{R}^2$. On the other hand, we cannot follow a curve for the slices of a specific variable $x_i$, if the subsystem of $n-1$ equations that we have chosen is independent of $x_i$ - that is, if the first $n-1$ elements of the $i$-th column of the Jacobian are zero.

Considering all the above, we have produced the additional algorithm:  $Reorder$, whose call can be inserted before the execution of $Step~3$ of Algorithm 2.1, and provides a suggestion for equation or variable reordering, that the user may follow. The routines $SwapRows$ and $SwapColumns$ can be used after the execution of $Reorder$, so that the final execution ordering lies on user decision.

{\vskip 0.3cm}
\noindent{\bf Algorithm 5.1: Reorder}\\{\it
Step 1: Input \{J\}\\
Step 2: n=length(J); D=zeros(n,n); rows=1:n; columns=1:n\\
Step 3: solvable=True; swappedrows=False; swappedcols=False\\
Step 4: Repeat Steps 5 - 7 for i=1,..,n\\
\hspace*{18pt} Step 5: Repeat Steps 6 - 7 for j=1,..,n\\
	\hspace*{35pt} Step 6: If Exists(x(j),J(i,j)): D(i,j)=2\\
	\hspace*{35pt} Step 7: Else If J(i,j)$\neq$ 0: D(i,j)=1\\
Step 8: Repeat Steps 9 - 11 for i=1,..,n-1\\
	\hspace*{18pt} Step 9: If sum(D(i,1:n-1))==0: Do Steps 10 - 11\\
	\hspace*{35pt} Step 10: If sum(D(n,1:n-1))==0 or swappedrows: solvable=False; End\\
	\hspace*{35pt} Step 11: Else: rows(n)=i; rows(i)=n; D=D(rows); swappedrows=True\\    
Step 12: If swappedrows: End\\
Step 13: For j=1,..,n CountOnes(j)=FindNumberOfOnes(D(1:n-1,j))\\
Step 14: MinOnes=min(CountOnes(j))\\ 
Step 15: Repeat Step 16 for j=1,..,n-1\\
\hspace*{18pt} Step 16: If CountOnes(j)==MinOnes: columns(j)=n; columns(n)=j; Exit\\  
Step 17: Output \{rows, columns, solvable\}
}
{\vskip 0.3cm} 
\noindent The basic idea behind steps 2-7 of Algorithm 5.1 is, to derive a matrix $D$ from $J$, whose $(i,j)$-element is:\\
\indent {\bf 2} if the function $f_i$ is nonlinearly dependent on variable $x_j$\\
\indent {\bf 1} if $f_i$ is linearly dependent on $x_j$, and\\ 
\indent {\bf 0} if $f_i$ is independent of $x_j$.
 
{\noindent} In steps 8-11, we eliminate the rank-reduction case of $J_l$, which may happen if a row of $D$ has all zeros but for the last element. This case can be solved by swapping the corresponding equation with the last equation of $F$ - and correspondingly swap this row with the last one in $D$ ($Step~11$). If another such row exists in $D$, then $rank(J)<n$ and therefore the system (1) cannot be uniquely solved and the flag variable $solvable$ is returned with the value $False$ ($Step~10$).

In fact, the occasion of zero elements everywhere but in the last element of a row in $D$, is the only occasion that Algorithm 5.1 suggests row-reordering. The reason behind this, is that all the linearity-in-a-variable cases can be treated by column swapping. After performing a great number of tests, we realised that row swapping byitself without column swapping did not improve significantly the method's performance. Nevertheless, we do provide row swapping as an additional option for the user to attempt.

In steps 13-16, we choose the first-encountered column ($j$) with the minimum number of ones in the rows $1..n-1$ and swap it with the last column of $D$. This corresponds to reordering the variables in (1) by swapping variables $x_j$ and $x_n$. We remind that the presence of 1 in an element of $D$ implies  linear dependence of the current function on the current variable. As noted earlier in this section, it is desirable such a variable to be included in the subsystem-to-follow, since it reduces the complexity of the problem. However, variable swapping can be done only once and only if row swapping has not happened before ($Step~12$). This is due to the fact that, after row swapping, a subsequent variable swapping may restore the zero-element case that row swapping has already cured.

The output $rows$ and $cols$ ordering of Algorithm 5.1 can be given as input to the routines $SwapRows$ and $SwapColumns$ respectively, in order to actually implement the algorithm's suggestion. In the next section we present results from applying our method to a variety of test-problems and illustrate various cases where the method's performance is improved by use of the ordering suggestions of Algorithm 5.1.

\section{Application to various test-problems with dimension up to $n=10$}
\label{sec6} 

In order to test our algorithm, we have applied it to more than 130 test-problems, in their majority taken from \cite{MGH}, \cite{Kear87} and \cite{Kuiken}. From those problems that are generalised for any dimension, we tested $n$ up to the value of 10. For all tested problems up to $n=10$, the algorithm managed to locate all the real solutions within a given $n$-box, with a wanted accuracy of $10^{-4}$ and with reasonable execution times (at most, of the order of seconds). As mentioned in Sect.~3, this accuracy is adequate for providing initial approximations of the solutions, which can subsequently be used as starting values by any high-precision and single-solution finding algorithm. 

\begin{table}[h]
\small
\centering
\caption{Application of our algorithm to the most representative test-problems up to $n=10$.}
\bigskip
\begin{tabular}{rrllrcl}
\hline
\#&$n$&reference&$n$-box&\#sols&secs&alg. sugg.\\
\hline
T1&10& Tridiagonal function \cite{Broyden}&$[-3,3]^{10}$&2&1.6&no reord.\\
T2&9&Almost-linear function \cite{Brown}&$[-20,20]^9$&3&0.2&no reord.\\
T3&8&R.K.P. \cite{Morgan}&$[-1,1]^8$&16&4.5&swap $x_5,x_8$\\
T4&7&D.I.E.F. \cite{More}&$[-5,5]^7$&1&0.1&no reord.\\
T5&6&EXP6 function \cite{Biggs}&$[-12,12]^6$&6&7.5&swap $x_1,x_6$\\
T6&5&Chebyquad function \cite{Fletcher}&$[0,1]^5$&120&165&no reord.\\
T7&4&System of quadratics \cite{Kear87}&$[-1,1]^4$&2&0.1&swap $x_1,x_4$\\
T8&3&Three-dimensional function \cite{Box}&$[0,11]^2\times[-2,2]$&84550&261&swap $x_1,x_3$\\
T9&2&Linear function \cite{MGH}&$[-5,5]^2$&1&0.04&swap rows\\
T10&2&Problem 1 \cite{Kuiken}&$[-1.6,1.6]\times[-1.04,1.04]$&12&0.9&no reord.\\
T11&2&Problem 2 \cite{Kuiken}&$[-3.1,3.8]\times[0.11,3.1]$&20&1.8&no reord.\\
\end{tabular}
\end{table}

In Table~1 we have chosen to present the application of our method to 11 test-problems, for values of $n$ from 10 down to 2, and trying to select, from our list of problems, either the most representative or the most troublesome. These problems we have renumerated, for the needs of later referring to them for extra commenting within the text. In all cases but for problem T8, the number of solutions presented in column $\#sols$ of Table~1, is the number of the problem's real solutions within the given $n$-box. This number coincides with the number of solutions found by our algorithm. Since the solutions of problem T8 are infinite, for this problem we indicatively present a number of solutions found by one execution of our method with specific input values (referred to in the text below). In column $secs$ we present the corresponding CPU-time needed for this execution. All experiments were performed on an Intel Core i5-3470, 3.60 GHz CPU, 8GB RAM, Windows 10 Pro x64 computer and all the codes were implemented in Matlab R2016a.

In what follows, we analyze each test-problem of Table~1.

{\bf T1}: A classical unconstrained optimization test-problem with variable $n$ is Broyden tridiagonal function \cite{Broyden} (see also test function 30 of \cite{MGH}). Here we present its function $F$ for $n=10$.
\begin{equation}
f_i=(3-2x_i)x_i-x_{i-1}-2x_{i+1}+1,~i=1,..,10,~x_0=x_{11}=0\\
\end{equation}
This problem has two real solutions in $[-3,3]^{10}$. The corresponding matrix $D$ of Algorithm 5.1 is tridiagonal, with 2 everywhere at the main diagonal and 1 everywhere at the lower and upper diagonals. Since it is the last column of $D$ that has the minimum number of ones, there is no need to do any variable swapping. The method is so well-behaved for this problem that it finds the two solutions even for step-sizes as large as $stepz=stepx=6$ (this is because it locates the solutions even from the first $z_0=-3$ slice). For $stepx=6$, $2^9=512$ mesh starting points are used and the CPU-time needed for the execution is $1.6''$. 

{\bf T2}: Another classical unconstrained optimization test-problem with variable $n$, is Brown almost-linear function \cite{Brown} (see also test function 27 of \cite{MGH}). Here we present its function $F$ for $n=9$. 
\begin{equation}
f_i=x_i+\sum_{j=1}^{9}x_j-10,~i=1,..,8,~f_9=\prod_{j=1}^{9}x_j-1
\end{equation}
All the solutions of this problem for any $n$, are of the form $(\alpha,...,\alpha,\alpha^{1-n})$ where $\alpha$ is any solution of the polynomial equation $na^n-(n+1)a^{n-1}+1=0$. For the particular case $n=9$, the equation $9a^9-10a^8+1=0$ has three real solutions: 
$\alpha=1$, $\alpha=0.974543355846$ and $\alpha=-0.7052133225$. These give rise to three corresponding real solutions of (6) and we can find them all if we search within the 9-box $[-20,20]^9$. The corresponding matrix $D$ of Algorithm 5.1 has all ones and therefore there is no need to do any row or variable swapping. Our method is very well-behaved also for this test-problem and can find the three solutions when using step-sizes as large as $stepz=stepx=40$, since it locates all the solutions even from the first $z_0=-20$ slice. For $stepx=40$, $2^8=256$ mesh starting points are used and the CPU-time needed for the execution is $0.2''$. 

{\bf T3}: This is a specific $n=8$ problem, mentioned as "Problem 11: A robot kinematics problem" (and referred to as R.K.P. in Table~1) by \cite{Kear87} and first used by \cite{Morgan}. Its function $F$ is
\begin{equation}
\begin{aligned}
f_1&=\alpha_1x_1x_3+\alpha_2x_2x_3+\alpha_3x_1+\alpha_4x_2+\alpha_5x_4+\alpha_6x_7+\alpha_7\\
f_2&=\alpha_8x_1x_3+\alpha_9x_2x_3+\alpha_{10}x_1+\alpha_{11}x_2+\alpha_{12}x_4+\alpha_{13}\\
f_3&=\alpha_{14}x_6x_8+\alpha_{15}x_1+\alpha_{16}x_2,~~f_4=\alpha_{17}x_1+\alpha_{18}x_2+\alpha_{19}\\
f_5&=x_1^2+x_2^2-1,~f_6=x_3^2+x_4^2-1,~f_7=x_5^2+x_6^2-1,~f_8=x_7^2+x_8^2-1\\
\end{aligned}
\end{equation}
where the values of the constants $\alpha_1,..,\alpha_{19}$ are given in \cite{Kear87}. The problem is known to have 16 real solutions in the 8-box $[-1,1]^8$. The corresponding matrix $D$ of Algorithm 5.1 is
\begin{equation*}
D=
\begin{bmatrix}
1&1&1&1&0&0&1&0\\
1&1&1&1&0&0&0&0\\
1&1&0&0&0&1&0&1\\
1&1&0&0&0&0&0&0\\
2&2&0&0&0&0&0&0\\
0&0&2&2&0&0&0&0\\
0&0&0&0&2&2&0&0\\
0&0&0&0&0&0&2&2\\
\end{bmatrix}
\end{equation*}
and the routine $Reorder$ suggests swapping of columns 5 and 8, i.e. of the variables $x_5$ and $x_8$ of the system, because column 5 has fewer ones than column 8. However, testing the system with or without the swapping of variables ends to the same performance. Actually, even without variable swapping, the method allows the use of very large step-sizes such as $stepz=stepx=2$. For $stepx=2$, $2^7=128$ mesh starting points are used and the method outputs all the 16 real solutions in a CPU-time of 4.5 seconds. Half of the solutions are found for the $z_0=-1$ slice, and the other half are found for the $z_0=1$ slice.

{\bf T4}: One more classical unconstrained optimization test-problem with variable $n$ (presented by \cite{MGH} as test function 29) is the "Discrete integral equation function" (referred to as D.I.E.F. in Table~1), first used as a test-problem by \cite{More}. Here we present its function $F$ for $n=7$. 
\begin{equation}
f_i=x_i+{1\over 16} \Bigl[(1-t_i)\sum_{j=1}^{i}t_j(x_j+t_j+1)^3+t_i\sum_{j=i+1}^{7}(1-t_j)(x_j+t_j+1)^3\Bigr]
,~t_i=i/8,~~i=1,..,7.
\end{equation}
In the 7-box $[-5,5]^7$ this problem has only one real solution.
Each element of the corresponding matrix $D$ of Algorithm 5.1 is 2, and so there is no need to do any row or variable swapping. Our method is very well-behaved also for this problem and the use of step-sizes can be as large as $stepz=stepx=10$, since it can locate the solution even from the first $z_0=-5$ slice. For $stepx=10$, $2^6=64$ mesh starting points are used and the CPU-time needed for the execution is $0.1''$. 

{\bf T5}: Another problem presented by \cite{MGH} (as test function 18) is Biggs EXP6 function \cite{Biggs} for $n=6$. The function is
\begin{equation}
\begin{aligned} 
f_i&=x_3e^{-t_ix_1}-x_4e^{-t_ix_2}+x_6e^{-t_ix_5}-y_i,\\
t_i&=0.1i,~~~y_i=e^{-t_i}-5e^{-10t_i}+3e^{-4t_i}, ~~~i=1,..,6.\\
\end{aligned}
\end{equation}
In the 6-box $[-12,12]^6$ this problem has 6 real solutions. Columns 1, 2 and 5 of the corresponding matrix $D$ of Algorithm 5.1 have all two's, while columns 3, 4 and 6 of $D$ have all ones. As a result, the routine $Reorder$ suggests swapping of columns 1 and 6. When testing our method without column swapping, it locates all 6 solutions with $stepz$ as large as 1.5 and $stepx$ as large as 6. The CPU time needed for the algorithm's execution with the above mentioned step-sizes is $25''$. Alternatively, if we perform the suggested swapping of columns 1 and 6 (that is, of variables $x_1$ and $x_6$), the method can locate all 6 solutions with $stepz$ as large as 3 and $stepx$ as large as 12, thus reducing the needed CPU-time to $7.5''$. As a result, test-problem T5 can be considered as an example of improvement of algorithm's performance through variable swapping. 

{\bf T6}: One more classical unconstrained optimization test-problem with variable $n$ presented by \cite{MGH} (as test function 35) is the Chebyquad function, first used as a test-problem in \cite{Fletcher}. The function $F$ is based on the shifted to $[0,1]$ Chebyshev polynomials $T_i(x)$. Here we present and use it for $n=5$. 
\begin{equation}
f_i={1\over 5} \sum_{j=1}^{5}T_i(x_j),~i=1,3,5,~~f_i={1\over 5} \sum_{j=1}^{5}T_i(x_j)+{1\over {i^2-1}},~i=2,4
\end{equation}
The above problem has $5!=120$ solutions in the 5-box $[0,1]^5$. The coordinates of these solutions in $\mathbb{R}^5$ are derived from all the different arrangements of the numbers 0.0838, 0.3127, 0.5000, 0.6873, 0.9162. The corresponding matrix $D$ of Algorithm 5.1 has ones in the first row and twos in all the other rows, therefore the algorithm does not suggest any row or variable swapping. The maximum distance between two $z_0$-slices within the $x_5$-interval $[0,1]$, so that all 24 curve-branches are located, is $stepz=0.005$. This necessarily small $stepz$ results to larger CPU-times compared to the other test-problems we have studied. E.g. for the combination of $stepz=0.005$ and $stepx=0.25$ (which creates a mesh of $5^4=625$ starting points) the method needs $165''$ to locate all the 120 solutions.  

{\bf T7}: Another variable-dimension problem is the system of quadratics presented in \cite{Kear87} (as test function 16) for $n=4$. The system's function is
\begin{equation}
f_i=(x_i-0.1)^2+x_{i+1}-0.1,~i=1,..,3,~f_4=(x_4-0.1)^2+x_1-0.1.
\end{equation}
In the 4-box $[-1,1]^4$ this problem has 2 real solutions. The corresponding matrix $D$ of Algorithm 5.1 is
\begin{equation*}
D=
\begin{bmatrix}
2&1&0&0\\
0&2&1&0\\
0&0&2&1\\
1&0&0&2\\
\end{bmatrix}
\end{equation*}
and the routine $Reorder$ suggests swapping of columns 1 and 4, i.e. of the variables $x_1$ and $x_4$ of the system, because column 1 has fewer ones than column 4 (if we discard the last row). Without column swapping, the method can approximate the solution with $x_i=0.1$ with accuracy only up to $10^{-2}$, despite using small step-sizes. On the contrary, if we perform the suggested swapping of variables $x_1$ and $x_4$, the method can calculate both solutions with accuracy as good as $10^{-14}$ and with step-sizes as large as $stepz=stepx=2$. For these step-sizes a mesh of $2^3=8$ starting points is used and the needed CPU-time is $0.1''$. Therefore, test-problem T7 is a profound example of improvement of algorithm's performance through variable swapping. 

{\bf T8}: A problem that requires special treatment is Box three-dimensional function \cite{Box}, also presented by \cite{MGH} as test function 12. The system's function is
\begin{equation}
f_i=e^{-t_ix_1}-e^{-t_ix_2}-x_3(e^{-t_i}-e^{-10t_i}),~~~t_i=0.1i,~~~i=1,2,3.
\end{equation}
This problem has the infinitely many real solutions (1,10,1),(10,1,-1),$(\alpha,\alpha,0)$, $\alpha\in\mathbb{R}$. 
We will investigate the performance of our method when searching for solutions in the 3-box $[0,11]^2\times[-2,2]$.

The corresponding matrix $D$ of Algorithm 5.1 has all twos in the first two columns and all ones in the third column. Consequently, the routine $Reorder$ suggests swapping of columns 1 and 3, i.e. of the variables $x_1$ and $x_3$ of the system, because column 1 has fewer ones than column 3. Furthermore, the nature of the problem's solutions is such, that the density of the solutions we can obtain, is affected by the density of the $x_n$-slices, only if we consider as "running", the variable $x_1$ or the variable $x_2$. As a result, when testing the algorithm without column swapping (that is, with $x_3$ as the running variable) and experimenting with the parameters $stepz$, $stepx$, $step$ and $thresh$ in order to locate as many solutions as possible (with $acc_2=10^{-10}$), we end up with at most 92 solutions. These are, the points $(1,10,1)$ and $(10,1,-1)$ plus 90 points of the form $(\alpha,\alpha,0)$. This is achieved in a CPU-time of $28''$ for the parameter values $stepz=1$, $stepx=step=threshx=0.1$. Decreasing any of these values further does not result in finding more solutions. 

On the other hand, when we implement the suggested by $Reorder$ routine column swapping, variable $x_1$ becomes the running variable. Then, the density of the located solutions is directly influenced by the magnitude of the $x_1$-slice and the $x_1$-step that is used in the $FollowCurve$ routine. After extensive testing it became obvious that, even if we sustain a big value of $stepz$, with adequately decreasing the values of the $step$ and $thresh$ curve-following parameters, we can regulate the density of the achieved solutions in the level we desire. For instance, starting with $stepz=stepx=1$ and $step=thresh=0.1$ (and for $acc_1=acc_2=10^{-10}$), we can locate the solution points $(1,10,1)$ and $(10,1,-1)$ plus 107 points of the form $(\alpha,\alpha,0)$ (with a density of order $O(10^{-1})$) in only $0.9''$. When decreasing  $step$ and $thresh$ to 0.001, we achieve the location of in total 8456 solutions with density of $O(10^{-3})$ in $14''$. Reducing the two last parameters even further to $step=thresh=10^{-4}$, the corresponding number of obtained solutions is increased to 84550. The CPU-time needed is now $261''$ and the solutions' density reaches the order $O(10^{-4})$. Therefore we can conclude that, for the test-problem T8, the need of column-swapping is unquestionable.

{\bf T9}: Another problem which, although trivial, requires special treatment is the Linear function of \cite{MGH}, whose equations for $n=2$ are
\begin{equation}
f_1=-x_2-1,~~f_2=-x_1-1.
\end{equation}
The corresponding matrix $D$ of Algorithm 5.1 is 
\begin{equation*}
D=
\begin{bmatrix}
0&1\\
1&0\\
\end{bmatrix}
\end{equation*}
and the routine $Reorder$ suggests swapping of rows 1 and 2, since otherwise the one-dimensional subsystem has zero-rank. As a result, if we attempt to execute the algorithm without row swapping and search for solutions e.g. in the 2-box $[-5,5]^2$, no solution is found. On the contrary, when we implement the suggested row swapping, the solution $(-1,-1)$ is found immediately from the first $x_2$-slice and in less than $0.1''$.

In what follows we will also study two more test-problems with $n=2$, the problems presented by Kuiken \cite{Kuiken}, so as to make a direct comparison with his 2D curve-following method.

{\bf T10}: This is the first problem of Kuiken \cite{Kuiken}, with equations: 
\begin{equation}
f_1=\Bigl(x_2-{1\over{3x_1}}\Bigr)\Bigl(x_2+{\tan}^{-1}(x_1)\Bigr),~~
f_2=\Bigl(x^2_2-{1\over{(1+x^2_1)^2}}\Bigr)\sin\Bigl({1\over{0.07+x^2_1+x^2_2}}\Bigr).
\end{equation}
This problem has 12 real solutions in the 2-box $[-1.6,1.6]\times[-1.04,1.04]$. The corresponding matrix $D$ of Algorithm 5.1 has all twos and therefore the algorithm does not suggest any row or variable swapping. Our method can locate all 12 solutions when using at most $stepz=stepx=0.7$ and the CPU time needed for this execution is only $0.9''$. On the contrary, for the same problem, Kuiken's method needs a much smaller step-size (0.02) for accurate curve-following, when using a second order Runge-Kutta scheme. This fact, together with the cost of coordinate system's rotations at each step of the curve-following, increases the CPU-time needed for Kuiken's method to obtain the same results to as much as $109''$.

{\bf T11}: This is the second problem of Kuiken \cite{Kuiken}, with equations: 
\begin{equation}
\begin{aligned} 
f_1&=\sin(1+x^2_1+x^2_2)-\cos(1+x^2_1+x^2_2){\tan}^{-1}(1+x^2_1+2x^2_2)\exp\Bigl({x^2_1+x^2_2\over 1+x^2_1+x^2_2}\Bigr)\\
f_2&=x^2_1\exp\Bigl({x^2_1-x^2_2\over 1+x^2_1+x^2_2}\Bigr)-\sqrt{\biggl|3x^2_1-2\exp\Bigl({x_1-x_2\over 1+|x_1|+|x_2|}\Bigr)\biggr|}.\\
\end{aligned} 
\end{equation}
This problem has 20 real solutions in the 2-box $[-3.1,3.8]\times[0.11,3.1]$. The corresponding matrix $D$ of Algorithm 5.1 has all twos and therefore the algorithm does not suggest any row or variable swapping. Our method can locate all 20 solutions when using at most $stepz=1.4,~stepx=0.6,~step=thresh=0.02$ and the CPU time needed for this execution is only $1.8''$. Again, for the same problem, Kuiken's method needs a step-size as small as 0.02 and the execution time reaches $427''$.

\section{Conclusions and suggestions for further work}
\label{sec7} 

In the present work we have presented a method for locating all real solutions of a system of nonlinear equations within a given $n$-box, by means of curve-following in $\mathbb{R}^n$. For the purposes of curve-following, Newton method has been used for numerically calculating a solution of the $(n-1)$-dimensional subsystem of the equations-to-follow, after starting from a previously found curve-point. The convergence of Newton's method is ensured by automatically adjusting the step-size. Newton's method was selected due to its quadratic convergence. However, any other solution-finding algorithm could be used in place of Newton and it would be interesting to investigate how this would affect the method's performance. In fact, even the present algorithm could be used recursively, for solving the $(n-1)$-dimensional subsystem. However, this would unnecessarily increase the computational cost, since what is needed at each curve-following step is the location of one solution in the proximity of the starting point, and not the location of all real solutions.

During the process of curve-following, the intersection point with the remaining hypersurface is found using bisection method. Again, here, any other more efficient solution-finding algorithm can be used instead, and we expect that this would affect the method's performance in a lesser or higher degree, depending on the number and density of the solutions to be found.

Various more sophisticated step-altering techniques might also be used to improve performance. In fact we have already tried to increase the step-size if $abs(u_1)>abs(u_0)$ at $Step~13$ of Algorithm 3.1 (which means that we are moving away from a solution). However, this technique proved not to  improve performance in the long term, since it deteriorated the accuracy of curve-following. In order to keep the accuracy within the user-specified limit $acc_1$, the step-size had to be halved again, thus leading to larger computational effort.   

Finally, a more complex row and column-ordering selection algorithm could also be developed, which a) could exploit the specific curve properties of each equation b) aim to enable the selection of the largest possible $stepx$ parameter. The reason for this, is that, the part of the algorithm 
with the highest complexity is the use of the $(n-1)$-dimensional starting points' mesh. As $n$ increases, it might be crucial for the method's performance, to select the row or variable-ordering that allows the use of the lowest-possible mesh density. This matter is still under investigation.

\section*{Acknowledgements}
The author would like to thank Prof. F. Petalidou for providing clarifications concerning geometrical notions.

\end{document}